\documentclass{amsart}
\usepackage{amssymb,latexsym,amsmath,amsthm}
\usepackage{fancyhdr}
  \address{School of Mathematics \\ Institute for Advanced Study \\ Einstein Drive, Princeton NJ 08540}
  \email{lbpierce@math.ias.edu}
\newtheorem{thm}{Theorem}
\newtheorem*{thm*}{Theorem}
\newtheorem{prop}{Proposition}
\newtheorem*{prop*}{Proposition}
\newtheorem{lemma}[prop]{Lemma}
\newtheorem*{lemma*}{Lemma}
\newtheorem{cor}{Corollary}[thm]
\newtheorem*{cor*}{Corollary}
\newtheorem*{conjecture*}{Conjecture}

\newtheorem{letterthm}{Theorem}

\theoremstyle{definition}

\newcommand{\bstack}[2]{\substack{#1 \\ #2}}

\newcommand{\maps}{\rightarrow}

\newcommand{\uhat}[1]{\check{#1}}

\newcommand{\lcm}{\text{lcm}}

\newcommand{\al}{\alpha}
\newcommand{\be}{\beta}

\newcommand{\ga}{\gamma}
\newcommand{\del}{\delta}
\newcommand{\ep}{\epsilon}

\newcommand{\Om}{\Omega}

\newcommand{\sig}{\sigma}

\newcommand{\Tcal}{\mathcal{T}}

\newcommand{\C}{\mathbb{C}}

\newcommand{\R}{\mathbb{R}}

\newcommand{\Z}{\mathbb{Z}}

\newcommand{\beq}{\begin{equation}}
\newcommand{\eeq}{\end{equation}}
\begin{document}

\title{A note on twisted discrete singular Radon transforms}
\author{Lillian B. Pierce}
\maketitle  

\begin{abstract}
In this paper we consider three types of discrete operators stemming from singular Radon transforms. We first extend an $\ell^p$ result for translation invariant discrete singular Radon transforms to a class of twisted operators including an additional oscillatory component, via a simple method of descent argument. Second, we note an $\ell^2$ bound for quasi-translation invariant discrete twisted Radon transforms. Finally, we extend an existing $\ell^2$ bound for a closely related non-translation invariant discrete oscillatory integral operator with singular kernel to an $\ell^p$ bound for all $1< p< \infty$. This requires an intricate induction argument involving layers of decompositions of the operator according to the Diophantine properties of the coefficients of its polynomial phase function.

\end{abstract}

\section{Introduction}

In this note we consider  twisted versions of discrete analogues of singular Radon transforms, and related discrete oscillatory integral operators. In their original setting, singular Radon transforms take the form
\beq\label{sRt}
 \Tcal f(x) = p.v. \int_{\R^{k_1}} f(\ga_t(x)) K(t) dt,
 \eeq
where $t\in \R^{k_1}$, $x \in \R^{k_2}$ and $\ga_t $ is a family of diffeomorphisms of $\R^{k_2}$ depending smoothly on $t$, such that $\ga_0$ is the identity. The kernel $K$ is a Calder\'{o}n-Zygmund kernel, that is, $K$ is $C^1$ on  $\R^{k_1} \setminus \{0 \}$ and satisfies the differential inequalities 
\beq\label{diff}
|\partial_t^\al K(t)| \leq A |t|^{-k_1 - |\al|}, \quad \text{for $0 \leq |\al | \leq 1$},
\eeq
and the cancellation condition
\beq\label{cancel}
| \int_{\ep \leq |t| \leq R} K(t)dt | \leq A
\eeq
uniformly in $0< \ep< R < \infty$. 

Such operators have been studied extensively, and the requirement on the varieties $\{ \ga_t(x) : t \in \R^{k_1} \}_{x \in \R^{k_2}}$ that guarantees the $L^p$ boundedness of $\Tcal$, namely a ``finite-type'' curvature condition with multiple equivalent formulations, is now understood (see \cite{CNSW} for the general theory).
Much less is known in the discrete setting, although significant progress has recently been made for two types of discrete singular Radon transforms: the translation invariant and quasi-translation invariant cases. 
\subsection{Translation invariant discrete operators}
Define for (compactly supported) functions $f : \Z^{k_2} \maps \C$ the discrete operator
\beq\label{T_P_dfn}
 T_Pf(n) = \sum_{\bstack{m \in \Z^{k_1}}{m \neq 0}} f(n-P(m))K(m),
 \eeq
where $K$ is a Calder\'{o}n-Zygmund kernel and $P=(P_1, \ldots, P_{k_2})$ is a polynomial mapping $\Z^{k_1} \maps \Z^{k_2}$. 
Such an operator was first considered by Arkhipov and Oskolkov \cite{AO}, who proved an $\ell^2$ result for dimension $k_1=k_2=1$, and by Stein and Wainger \cite{SW0}, who proved that $T_P$ is bounded on $\ell^p$ for $3/2< p < 3,$ in all dimensions. A recent deep result of Ionescu and Wainger \cite{IW} proves the desired $\ell^p$ bounds for $T_P$ for all $1<p< \infty$ and all dimensions: 
\begin{letterthm}\label{IW_thm}
The operator $T_P$ extends to a bounded operator on $\ell^p(\Z^{k_2})$ for $1< p < \infty$, with
\[ ||T_Pf||_{\ell^p(\Z^{k_2})} \leq A_p ||f||_{\ell^p(\Z^{k_2})}.\]
The constant $A_p$ depends only on $p$, the dimension $k_1$, and the degree of the polynomial $P$.
\end{letterthm}

The first result of this note is an extension of Theorem \ref{IW_thm} to the class of twisted translation invariant discrete singular Radon transforms. Let $T_{P,Q}$ be the operator
\beq\label{TPQ}
T_{P,Q}f(n) = \sum_{\bstack{m \in \Z^{k_1}}{m \neq 0}} f(n-P(m)) K(m) e^{2\pi i Q(m)},
\eeq
where $P$ is again a polynomial mapping $\Z^{k_1} \maps \Z^{k_2}$, $Q$ is a polynomial mapping $\R^{k_1} \maps \R$, and $K$ is a Calder\'{o}n-Zygmund kernel.
We prove: 
\begin{thm}\label{IW_gen}
The operator $T_{P,Q}$ extends to a bounded operator on $\ell^p (\Z^{k_2})$ for $1< p< \infty$, with
\[ ||T_{P,Q}f||_{\ell^p(\Z^{k_2})} \leq A_p ||f||_{\ell^p(\Z^{k_2})}.\]
The constant $A_p$ depends only on $p$, the dimension $k_1$, and the degrees of $P$ and $Q$.
\end{thm}

Note that in both the above theorems, the operator norm is independent of the coefficients of the polynomials $P$ and $Q$; this is a feature of all the results in this paper. The proof of Theorem \ref{IW_thm} due to Ionescu and Wainger is technically impressive, involving ``almost orthogonality'' properties of an intricate decomposition of the Fourier multiplier of the operator $T_P$ according to Diophantine approximations of the spectral variables, motivated by ideas from the circle method of Hardy and Littlewood. Indeed, many of the current approaches to discrete operators are rooted in circle method techniques: for example, the pioneering work of Bourgain \cite{Bour88A} \cite{Bour89}, further developed by Ionescu, Magyar, Stein and Wainger in \cite{SW0} \cite{SW1} \cite{SW2} \cite{SW3} \cite{MSW} \cite{IW} \cite{IMSW}. The surprising aspect of Theorem \ref{IW_gen} is that we  may ultimately avoid the beautiful but substantial technical work of the original method, instead proving Theorem \ref{IW_gen} via a transference principle or ``method of descent,'' which allows us to inject the twist into the Ionescu and Wainger paradigm quite simply.

\subsection{Quasi-translation invariant discrete operators}\label{sec_IMSW}
 Non-translation invariant singular Radon transforms of the form (\ref{sRt}) present new difficulties, even in the continuous setting, as the lack of translation invariance severely hampers the use of the Fourier transform. One may easily define a discrete analogue of such operators, but the inapplicability of Fourier transform methods seems particularly problematic in the discrete setting: circle method techniques rely on a decomposition of spectral variables.
 
 Thus we restrict our attention to the more tractable ``quasi-translation invariant'' operators acting on (compactly supported) functions $f : \Z^{k}\times \Z^l \maps \C$ by
\[R_Pf(n, n') = \sum_{\bstack{m \in \Z^k}{m\neq 0}} f(n-m, n'-P(n,m))K(m),\]
where $(n, n') \in \Z^k \times \Z^l$, $m\in \Z^k$,  $P$ is a polynomial mapping $\Z^k \times \Z^k \rightarrow \Z^l$, and $K$ is a Calder\'{o}n-Zygmund kernel. 
In \cite{SW1} Stein and Wainger proved an $\ell^2$ bound for operators of the form $R_P$, which has recently been extended in \cite{IMSW} to $\ell^p$ bounds for all $1<p< \infty$, albeit with a restriction on the degree of $P$:
\begin{letterthm}\label{SW1_thm}
If the polynomial $P$ is of degree at most 2, the operator $R_P$ extends to a bounded operator on $\ell^p (\Z^k \times \Z^l)$ for $1< p < \infty$, with
\[ ||R_Pf||_{\ell^p(\Z^k \times \Z^l)} \leq A_p ||f||_{\ell^p(\Z^{k} \times \Z^l)} .\]
The constant $A_p$ depends only on $p$ and the dimension $k.$  For $p=2$, this result holds for polynomials $P$ of any degree, in which case the constant $A_2$ depends on the dimension $k$ and the degree of the polynomial $P$. 
\end{letterthm}

The study of the operator $R_P$ on $\ell^2$ leads to a related discrete oscillatory integral operator:
\beq\label{SW_I_dfn}
 Tf(n) = \sum_{\bstack{m \in \Z^k}{n-m \neq 0}} f(m) K(n-m)e^{2\pi i Q(n,m)},
 \eeq
where $n \in \Z^k$, $Q$ is a polynomial mapping $\R^k \times \R^k \maps \R$, and $K$ is a Calder\'{o}n-Zygmund kernel. The main result of \cite{SW1} is that the operator $T$ is bounded on $\ell^2 (\Z^k)$, with a bound dependent only on the degree of $Q$, and independent of its coefficients; by Plancherel's theorem, this provides an $\ell^2$ bound for the operator $R_P$.

In this note we extend the $\ell^2$ result for $T$ to an $\ell^p$ result for all $1< p< \infty$; this may be seen as a discrete analogue of the results of Ricci and Stein \cite{RSI} for oscillatory integral operators on $\R^k$.
\begin{thm}\label{SW_gen}
The operator $T$ extends to a bounded operator on $\ell^p (\Z^k)$ for $1< p< \infty$, with
\[ ||Tf||_{\ell^p(\Z^k)} \leq A_p ||f||_{\ell^p(\Z^{k})} .\]
The constant $A_p$ depends only on $p$, the dimension $k$, and  the degree of $Q$.
\end{thm}

To prove this, we follow the method of proof of \cite{SW1}, utilizing an inductive decomposition of the operator $T$ based on the Diophantine properties of the coefficients of the polynomial $Q$, ultimately reducing the operator to a tensor product of a discrete Gauss sum operator with nice arithmetic properties, and a discrete operator that may be closely approximated by its continuous analogue, a singular integral operator with oscillatory kernel, which is known to be bounded on $L^p$ by  \cite{RSI}. 

The third result of this note is the observation that for $p=2$, Theorem \ref{SW1_thm}  extends immediately to the twisted case.
Let $R_{P,Q}$ be the operator
\beq\label{R_tilde_op}
 R_{P,Q} f(n, n') = \sum_{\bstack{m \in \Z^k}{m\neq 0}} f(n-m, n'-P(n,m))K(m)e^{2\pi i Q(n,m)},
 \eeq
where $(n, n') \in \Z^k \times \Z^l$, $m\in \Z^k$, $P$ is a polynomial mapping $\Z^k \times \Z^k \rightarrow \Z^l$, $Q$ is a polynomial mapping $\R^k \times \R^k \rightarrow \R$, and $K$ is a Calder\'{o}n-Zygmund kernel.
\begin{cor}\label{SW_cor}
The operator  $R_{P,Q}$ extends to a bounded operator on $\ell^2 (\Z^k\times \Z^l)$, with
\[||R_{P,Q}f||_{\ell^2(\Z^k \times \Z^l)} \leq A ||f||_{\ell^2(\Z^{k} \times \Z^l)} .\]
The constant $A$ depends only on the dimension $k$ and the degrees of $P,Q$.
\end{cor}

In Section \ref{sec_IW_gen} we prove Theorem \ref{IW_gen} as a consequence of Theorem \ref{IW_thm}. In Section \ref{sec_SW_cor} we briefly note the proof of Corollary \ref{SW_cor}, before then turning in Section \ref{sec_SW_gen} to the proof of Theorem \ref{SW_gen}, which is the focus of the remainder of the paper. In Section \ref{sec_induction} we detail the inductive procedure that allows us to reduce the theorem to bounding a Gauss sum operator and a discrete oscillatory singular integral operator, which we then do in Section \ref{sec_final_bd}, thus completing the proof of the theorem.

In what follows, the discrete Fourier transform of a function $f \in \ell^1(\Z^{k})$ is defined by
\[ \hat{f} (\xi) = \sum_{n \in \Z^k} f(n) e^{-2\pi i n \cdot \xi},\]
where $\xi \in \R^k$. Since $\hat{f} (\xi) = \hat{f}(\xi + n)$ for any $n \in \Z^k$, it is a periodic function. For periodic functions $h \in L^2_{loc} (\R^k)$, we define the Fourier inverse to be
\[ \uhat{h} (n) = \int_{[0,1]^k} h(\xi) e^{2\pi i \xi \cdot n} d\xi.\]

%*********************************************
\section{Proof of Theorem \ref{IW_gen}}\label{sec_IW_gen}
%*********************************************
We first consider the translation invariant operator $T_P$ given by (\ref{T_P_dfn}).  
In fact, it is equivalent to bound $T_P$ on $\ell^p(\Z^{k_2})$ as to prove an $L^p(\R^{k_2})$ bound for the related operator acting on functions of $\R^{k_2}$ by 
\[ T^\sharp_PF(x) = \sum_{\bstack{m \in \Z^{k_1}}{m \neq 0}} F(x-P(m))K(m).\]
A similar equivalence holds for $T_{P,Q}$ defined by (\ref{TPQ}) and the operator 
\[T^\sharp_{P,Q}F(x) = \sum_{\bstack{m \in \Z^{k_1}}{m \neq 0}} F(x-P(m))K(m)e^{2\pi i Q(m)}.\]
This is simply a consequence of the fact that $P(\Z^{k_1}) \subseteq \Z^{k_2}$, and may be seen, for example in the case of $T^\sharp_P$, as follows. (Similar arguments, stemming from an observation of E. M. Stein, arise in \cite{IW} and \cite{IMSW}.) 
Given a function $f$ defined on $\Z^{k_2}$, define $F$ on $\R^{k_2}$ by setting $F(x) = f(n)$ for $x$ in the unit cube centered at $n \in \Z^{k_2}$. Precisely, setting $Q = (-1/2,1/2]^{k_2}$ and letting $\chi_Q$ denote the characteristic function of $Q$,  $F$ is defined by
\[F(x) = \sum_{n \in \Z^{k_2}} f(n) \chi_{Q}(x-n).\]
Note that $F$ belongs to $L^p(\R^{k_2})$ precisely when $f \in \ell^p(\Z^{k_2})$, and moreover  that $||F||_{L^p(\R^{k_2})} = ||f||_{\ell^p(\Z^{k_2})}$. Furthermore, it is immediate that $T_P^\sharp F(x) = T_Pf(n)$ for $x \in Q+n$, so that
\[ ||T_Pf ||_{\ell^p(\Z^{k_2})}^p =  \sum_n | T_Pf(n)|^p
	 = \sum_n \int_{Q+n} | T^\sharp_PF(x)|^p dx
	 	= ||T_P^\sharp F||_{L^p(\R^{k_2})}^p,\]
and hence $||T_P||_{\ell^p} \leq ||T_P^\sharp ||_{L^p}$.

Conversely, given a function $F$ on $\R^{k_2}$, define for each $x \in Q $ a function $f_x$ acting on $\Z^{k_2}$ by $f_x(m) = F(x+m)$ for all $m \in \Z^{k_2}$. 
Then
\begin{eqnarray*}
||T_P^\sharp F||_{L^p}^p 	& = &
% \int_{\R^k} | \sum_m K(m) F(x-P(m)) |^p dx 
	 	 \sum_n \int_Q | \sum_m K(m) F(x+n-P(m) )|^p dx 
		%&= & \sum_n \int_Q | \sum_m K(m) f_x(n-P(m))|^p dx 
			= \sum_n \int_Q |T_P f_x(n)|^p dx \\
		& = &
			% \int_Q \sum_n |T_P f_x(n)|^p dx 
			  \int_Q ||T_P f_x(\cdot)||_{\ell^p}^p dx
			 \leq ||T_P||^p_{\ell^p} \int_Q ||f_x(\cdot)||_{\ell^p}^p  dx
			 = ||T_P||^p_{\ell^p}  \int_Q \sum_n |f_x(n)|^p dx
			%= ||T_P||^p_{\ell^p} \sum_n \int_Q |f_x(n)|^p dx
			 \\
			&=& ||T_P||^p_{\ell^p} \sum_n \int_Q |F(x+n)|^p dx 
		%= ||T_P||^p_{\ell^p}  \int_{\R^k} |F(x)|^p dx 
		=||T_P||^p_{\ell^p}\;  ||F||_{L^p}^p.		
\end{eqnarray*}
Therefore $||T_P^\sharp||_{L^p} \leq ||T_P||_{\ell^p}$ and the equivalence follows; a similar argument holds for $T_{P,Q}^\sharp$ and $T_{P,Q}$.

In conclusion, in order to prove Theorem \ref{IW_gen}, it suffices to consider the operators $T^\sharp_P$ and $T^\sharp_{P,Q}$. Applying the Euclidean Fourier transform, the operator $T^\sharp_P$ satisfies $(T^\sharp_P f)\hat{\;} (\xi)  = m(\xi) \hat{f}(\xi)$, with Fourier multiplier 
\[m(\xi) = \sum_{\bstack{m \in \Z^{k_1}}{m \neq0}} K(m) e(- \xi \cdot P(m)),\]
where $\xi \in \R^{k_2}$ and $e(t)$ denotes $e^{2\pi i t}$. 
The twisted operator $T_{P,Q}^\sharp$ has Fourier multiplier
\beq\label{twisted_mult}
\tilde{m}(\xi) = \sum_{\bstack{m \in \Z^{k_1}}{m \neq0}} K(m) e(- ( \xi \cdot P(m) - Q(m))),
\eeq
where again $\xi \in \R^{k_2}$.
The key to proving Theorem \ref{IW_gen} is showing that both $m(\xi)$ and $\tilde{m}(\xi)$ may be reduced to a ``universal'' multiplier via a transference principle, which we now record (as proved in Chapter 11 \S 4.6 of \cite{SteinHA}):
\begin{lemma}[Method of descent]\label{descent}
Let $L: \R^{n_1} \maps \R^{n_2}$ be a linear map and $m :\R^{n_2} \maps \C$ a continuous function. Define the function $m_L : \R^{n_1} \maps \C$ by $m_L(\xi) = m(L\xi)$. Then for any $1<p<\infty$, the $L^p(\R^{n_1})$ norm of the operator defined by the Fourier multiplier $m_L$ does not exceed the $L^p(\R^{n_2})$ norm of the operator defined by the Fourier multiplier $m$:
\[ ||m_L||_{\mathcal{M}_p(\R^{n_1})} \leq ||m||_{\mathcal{M}_p(\R^{n_2})}.\]
\end{lemma}

In their original argument, Ionescu and Wainger employ this principle to reduce the given polynomial $P: \Z^{k_1} \maps \Z^{k_2}$ present in the multiplier $m$ to a ``generic polynomial'' as follows. Suppose that $P=(P_1,\ldots, P_{k_2})$ is of degree $d_P$ and each component $P_l$ is given by 
\[P_l(x) = \sum_{1 \leq |\al| \leq d_P} \beta_{l,\al}x^\al.\]
 (Note that we may disregard terms of order zero.) Set $D_P$ to be the cardinality of the set of multi-indices of order no more than $d_P$, 
 \[ D_P = |\text{Ind}(d_P)| = |\{\al \in \Z^{k_1}_{\geq 0} : 1 \leq |\al| \leq d_P\}|,\]
 where the norm of a multi-index is given by $|\al| = \al_1 + \cdots + \al_{k_1}$.
Denote $\xi \in \R^{D_P}$ componentwise by $(\xi_\al) = ([\xi]_\al)$, where $\al \in \text{Ind}(d_P)$ and $\text{Ind}(d_P)$ is ordered in a fixed manner. 
Then, define the \emph{universal multiplier} $\mu : \R^{D_P} \maps \C$ by
\beq\label{univ_mult_dfn}
 \mu(\xi) =\sum_{\bstack{m \in \Z^{k_1}}{m \neq0}}K(m)e(- \sum_{1 \leq |\al| \leq {d_P}} m^\al \xi_\al ).
 \eeq
To prove an \emph{a priori} inequality, one may assume without loss of generality that the kernel $K$ is compactly supported and thus the multiplier $\mu$ is continuous. 

Define a linear map $L : \R^{k_2} \maps \R^{D_P}$ componentwise by
\beq\label{lin_map}
 [L \eta]_\al = \sum_{l=1}^{k_2} \beta_{l,\al} \eta_l.
 \eeq
It is simple to check that the multiplier $\mu_L: \R^{k_2} \maps \C$ defined by $\mu_L(\xi) = \mu(L\xi)$, with $\mu$ the universal multiplier, and $L$ the linear map defined by (\ref{lin_map}), is precisely the Fourier multiplier $m$ of the operator $T_P^\sharp$. Thus by Lemma \ref{descent},
\[  ||m||_{\mathcal{M}_p(\R^{k_2})} =||\mu_L||_{\mathcal{M}_p(\R^{k_2})} \leq ||\mu||_{\mathcal{M}_p(\R^{D_P})}.\]
As a result, in order to bound the operator $T_P^\sharp$ on $L^p(\R^{k_2})$, it is sufficient to bound the operator with universal multiplier $\mu$ on $L^p(\R^{D_P})$. 
 Note furthermore that the definition (\ref{univ_mult_dfn}) of $\mu$ could also be written as 
\beq\label{mult_dot}
 \mu(\xi)  = \sum_{\bstack{m \in \Z^{k_1}}{m \neq0}}K(m) e(-\xi \cdot P_0(m)),
 \eeq
where $\xi \in \R^{D_P}$ and  $P_0 : \R^{k_1} \maps \R^{D_P}$  is the \emph{generic polynomial of degree $d_P$} with all coefficients equal to 1, defined componentwise by $[P_0(x)]_\al = x^\al$.
The problem has thus been reduced to considering the operator $T_{P_0}^\sharp$ defined by (\ref{T_P_dfn}) with $P_0$ in place of $P$; this Ionescu and Wainger do with great finesse, proving that $T_{P_0}^\sharp$ is bounded on $L^p(\R^{D_P})$ for all $1 < p < \infty.$ 

Note that this process increases, possibly significantly, the dimension of the underlying space, but we have gained an advantage by eliminating the role of the specific coefficients of $P$. 
We will now further use the method of descent to insert the extra oscillatory component $e^{2\pi i Q(m)}$ present in the operator $T_{P,Q}^\sharp$. Let $d_Q$ be the degree of the polynomial $Q$ and $d_P$ the degree of $P$. Set $d = \max (d_P,d_Q)$ and set $D = |\text{Ind}(d)| = |\{\al \in \Z^{k_1}_{\geq 0} : 1 \leq |\al| \leq d \}|$. It is convenient to use the convention that if $d_Q>d_P$, we include terms up to degree $d_Q$ in $P$ with zero coefficients, and vice versa; thus from now on we will consider both $P,Q$ to have degree $d$ in this sense. For future reference, note that we could repeat the above procedure for the original multiplier $m$ after padding $P$ with zeroes to bring the degree of $P$ up to $d$; this would replace each instance of $d_P$ by $d$ and each instance of $D_P$ by $D$. In conclusion, we would again reduce the problem to considering the universal multiplier $\mu$ in (\ref{mult_dot}), but now with dimension $D$ in place of $D_P$, and the generic polynomial $P_0$ of degree $d$.

We now define the \emph{quasi-universal multiplier} $\tilde{\mu}(\xi)$ for $\xi \in \R^D$ by
\[ \tilde{\mu}(\xi) =\sum_{\bstack{m \in \Z^{k_1}}{m \neq0}} K(m) e(- \sum_{1 \leq |\al| \leq d} m^\al \xi_\al )e( Q(m)) .\]
While this multiplier has replaced $P$ by the generic polynomial, it retains the specific polynomial $Q$. Applying the linear operator $L$ defined  in terms of the coefficients of $P$ by (\ref{lin_map}), 
we then note that as multipliers, 
\[\tilde{\mu}_L(\xi) = \tilde{\mu}(L\xi) = \tilde{m}(\xi),\]
where $\tilde{m}$ is the multiplier (\ref{twisted_mult}) of the twisted operator $T_{P,Q}^\sharp$.
 Furthermore, if $Q$ is defined by $Q(x) = \sum_{1 \leq |\al| \leq d} \theta_\al x^\al$, define the vector $\theta \in \R^{D}$ of coefficients of $Q$ componentwise by $[\theta]_\al = \theta_\al$. 
  Then by the definition of the generic polynomial $P_0: \R^{k_1} \maps \R^D$, $Q(x) = \theta \cdot P_0(x)$. Thus in fact our quasi-universal multiplier $\tilde{\mu}$ may be written as 
\beq\label{tilde_mult_dot} \tilde{\mu}(\xi) =\sum_{\bstack{m \in \Z^{k_1}}{m \neq0}} K(m) e(-(\xi - \theta)\cdot P_0(m) ) .
\eeq
Finally, by Lemma \ref{descent}, 
\[||\tilde{m}||_{\mathcal{M}_p(\R^{k_2})} = ||\tilde{\mu}_L||_{\mathcal{M}_p(\R^{k_2})} \leq ||\tilde{\mu}||_{\mathcal{M}_p(\R^D)}.\]
Thus in order to bound $T_{P,Q}^\sharp$ on $L^p(\R^{k_2})$ it is sufficient to bound the operator with Fourier multiplier $\tilde{\mu}$ on $L^p(\R^D)$.

At this point, it is possible to use the method of Ionescu and Wainger \cite{IW}, applied to $\tilde{\mu}$, to prove Theorem \ref{IW_gen} directly. But far simpler, 
comparing (\ref{mult_dot}) to (\ref{tilde_mult_dot}), it is clear that 
\[ \tilde{\mu} (\xi) = \mu(\xi-\theta). \]
It is helpful to keep in mind that the vector $\theta$ is fixed, once and for all, by the coefficients of the polynomial $Q$, while $\xi$ is the spectral variable. 
Recall that the operator with Fourier multiplier $\mu(\xi)$ is $T_{P_0}^\sharp$, so that we may write 
\[T_{P_0}^\sharp f (x) = (f * K_0) (x),\] 
where $K_0 (x) = \uhat{\mu} (x)$. Then the operator with multiplier $\mu(\xi-\theta)$, with fixed $\theta$, has kernel $K_0(x) e^{2\pi i \theta \cdot x}$; 
thus we have reduced bounding the operator $T_{P,Q}^\sharp$ to bounding the operator 
\[\tilde{T}_{P_0}^\sharp f(x)  = (f * K_0(m)e^{2\pi i \theta \cdot m}) (x).\]
But in fact $||\tilde{T}_{P_0}^\sharp f||_{L^p} = ||T_{P_0}^\sharp f||_{L^p}$ for all $1< p < \infty$, since by definition
\[ \tilde{T}_{P_0}^\sharp f(x)  = \sum_{m} f(m) K_0(x-m) e^{2\pi i (x-m) \cdot \theta},\]
and oscillatory factors that depend only on $m$ or only on $x$ do not affect the norm of the operator. 
Thus the $L^p$ boundedness of $\tilde{T}_{P_0}^\sharp$, and hence of $T_{P,Q}^\sharp$, follows immediately from Theorem \ref{IW_thm} of Ionescu and Wainger. This concludes the proof of Theorem \ref{IW_gen}.

%********************************************
\section{Proof of Corollary \ref{SW_cor}}\label{sec_SW_cor}
%********************************************

We now turn to quasi-translation invariant Radon transforms. 
We first record the simple observation that Corollary \ref{SW_cor} follows immediately from Theorem \ref{SW_gen}.
Taking the discrete Fourier transform of the operator $R_{P,Q}$ defined in (\ref{R_tilde_op}) with respect to the variable $n'$ alone gives
\[ (R_{P,Q} f)^{\hat{n'}}(n,\xi) =  \sum_{\bstack{m \in \Z^k}{m \neq 0}} e^{2\pi i (Q(n,n-m) - \xi \cdot P(n,n-m))} K(n-m)f^{\hat{n'}}(m,\xi),\]
where $\xi \in [0,1]^{l}$.
This operator is now of the form (\ref{SW_I_dfn}), with polynomial $Q_0(n,m) = Q(n,n-m) - \xi \cdot P(n,n-m)$. Noting that the operator norm provided by Theorem \ref{SW_gen} is independent of the coefficients of the polynomial $Q_0$, and hence of $\xi$, Corollary \ref{SW_cor} then follows by Plancherel's theorem in the second variable. 

We note that for $p \neq 2$, extending the $\ell^p$ result of \cite{IMSW} recorded in Theorem \ref{SW1_thm} to the twisted operator $R_{P,Q}$ is complicated by the manner in which the existing result is proved, involving passing through a Lie group defined in terms of the specific polynomial $P$. This method currently appears to admit a generalization to $\ell^p$ bounds for $R_{P,Q}$ only for certain pairs of polynomials $P, Q$, and thus we do not present it here.

%********************************************
\section{Proof of Theorem \ref{SW_gen}}\label{sec_SW_gen}
%********************************************

We now turn to Theorem \ref{SW_gen} and the operator $T$ defined in (\ref{SW_I_dfn}). For this we employ a double decomposition, first with respect to a dyadic decomposition of the Calder\'{o}n-Zygmund kernel, and then with respect to the Diophantine properties of the coefficients of the polynomial $Q$, following the presentation of \cite{SW1}. Whereas  $T$ is  not bounded on $\ell^\infty$, the multitude of operators we will encounter in the decomposition will be bounded on $\ell^\infty$, enabling us to interpolate with existing nontrivial $\ell^2$ bounds in order to capture $\ell^p$ results. As the original argument of \cite{SW1} is quite elaborate, we will focus only on the components that are necessary for Theorem \ref{SW_gen}, and merely sketch components of the proof not affected by the $\ell^p$ context.

Recall that 
\beq\label{T_dfn}
 Tf(n) = \sum_{\bstack{m \in \Z^k} {n-m \neq 0}} e^{2\pi i Q(n,m)}K(n-m) f(m).
 \eeq

The properties (\ref{diff}) and (\ref{cancel}) of the kernel $K$ allow its decomposition for $|x| \geq 1$ (the only region of interest for the discrete problem) into
\[ K(x) = \sum_{j=0}^\infty K_j(x),\]
where each $K_j$ is supported in $2^{j-1} \leq |x| \leq 2^{j+1}$, satisfies the same differential inequalities (\ref{diff}) as $K$, uniformly in $j$, and satisfies the mean value property $\int K_j(x)dx =0$. (For a derivation of this decomposition, see for example Chapter 6 \S 4.5 and Chapter 7 \S3.4 of \cite{SteinHA}.)

Similarly, we decompose the operator as  $T = \sum_{j=0}^\infty T_j,$
where 
\[T_jf(n) = \sum_{\bstack{m \in \Z^k}{n-m \neq 0}} e^{2\pi i Q(n,m)}K_j(n-m) f(m).\]
We now proceed with an inductive decomposition that classifies the indices $j$ in terms of the Diophantine properties of the coefficients of the polynomial $Q$, beginning with the highest degree terms.

\subsection{Outline of the inductive major/minor decomposition}
Suppose that 
\[Q(n,m) = \sum \theta_{\al, \be}n^\al m^\be,\]
 where the multi-indices $\al,\be$ satisfy $2 \leq |\al| + |\be| \leq d$, $d$ being the degree of the polynomial. Note that we may assume in each term that both $|\al| \neq 0$ and $|\be| \neq 0$: any oscillatory factor purely in terms of $n$ may be pulled out of the sum (\ref{T_dfn}) without affecting the norm of the operator, while any oscillatory factor purely in terms of $m$ may be cancelled by pre-multiplication of the function by an exponential factor. Since $Q(n,m)$ is a phase function, we may also assume that $\theta_{\al,\be} \in (0,1)$ for all $\al, \be$.

The decomposition begins with the coefficients $\theta_{\al,\be}$ of highest degree, with $|\al| + |\be| = d$. Fix $j \geq 0$, and a small number $\ep_d>0$. By the Dirichlet approximation principle, for each such $\theta_{\al,\be}$ there exist integers $ a_{\al,\be}=a_{\al,\be, j}$ and $ q_{\al,\be}=q_{\al,\be, j} $ such that
\beq\label{theta_aqj}
 \left| \theta_{\al,\be} - \frac{a_{\al,\be}}{q_{\al,\be}} \right| \leq \frac{1}{q_{\al,\be}2^{(d-\ep_d)j}},
 \eeq
where $1 \leq q_{\al,\be} \leq 2^{(d-\ep_d)j}$ and $1 \leq a_{\al,\be} \leq q_{\al, \be}$ with $(a_{\al,\be}, q_{\al,\be})=1$. We will call the set of all $q_{\al, \be}$ chosen in this manner, for all $|\al| + |\be| =d$ (and fixed $j$), the set of \emph{denominators of level $d$}. 

We now distinguish between two cases: major and minor indices $j$. In the first case, at least one of the denominators of level $d$ is ``large,'' specifically $q_{\al,\be} >2^{\ep_dj}$, in which case we call $j$ a \emph{minor index}. In the second case, all of the denominators of level $d$ are ``small,'' specifically $1 \leq q_{\al,\be}  \leq 2^{\ep_dj}$, and we call $j$ a \emph{major index}. To be more specific, we can also say that $j$ is \emph{major/minor of level $d$.}\footnote{Note that in a more typical circle method decomposition of an operator, for each fixed $j$ the spectral variable would be classified as belonging to a major or minor arc with respect to $j$. However, in this case we do not consider a Fourier multiplier and $\theta$ is not a variable but a fixed vector of real coefficients, and so we proceed in the opposite direction and decompose $j$ with respect to $\theta$.}

We then decompose the operator as 
\[ T = T_{M} + T_{m},\]
where $T_M$ is the sum of $T_j$ for all major $j$ of level $d$, and $T_m$ is the sum of $T_j$ for all minor $j$ of level $d$. 
The next step is to give a nontrivial $\ell^p$ estimate for each $T_j$ with $j$ minor, of the form
\[ ||T_j||_{\ell^p(\Z^k)} \leq A2^{-\del j} \]
for some $\del >0$. (For notational convenience we will write $||T||_{\ell^p(\cdot)}$ for $||T||_{\ell^p(\cdot) \maps \ell^p(\cdot)}$, for norms of operators that preserve an $\ell^p$ space.) Such an estimate would immediately allow one to sum over all minor $j$ to obtain a bound for the full minor component of the form
\[ ||T_m||_{\ell^p(\Z^k)} \leq A.\]

The main analysis goes into evaluating the major component $T_M$ via an inductive procedure, which at each stage breaks the existing major component at level $l$ into sub-components, the major and minor parts of level $l-1$, based on the Diophantine properties of the coefficients of terms of degree $l-1$ in the polynomial $Q$. We begin by outlining the transition from level $d$ to level $d-1$. First, write $T_M$ as 
\beq\label{TaqM}
T_M = \sum_{(a/q)_d} T^{(a/q)_d},
\eeq
 where we define 
\beq\label{Taq}
 T^{(a/q)_d} = \sum_{j \; \text{major}} T_j.
 \eeq
Here each $(a/q)_d$ denotes a specific collection $\{a_{\al,\be},q_{\al,\be}\}$ of pairs of rationals, where $\al, \be$ range over $|\al| + |\be|=d$ and the summation (\ref{Taq}) is restricted to those $j$ 
for which $a_{\al,\be, j} = a_{\al,\be}$, $q_{\al,\be,j}=q_{\al,\be}$ in the approximation (\ref{theta_aqj}), for all $|\al| + |\be| = d$, and such that
all the denominators satisfy $q_{\al, \be} \leq 2^{\ep_d j}$. The sum over $(a/q)_d$ in (\ref{TaqM}) then ranges over all such collections of rationals. The goal is then to bound the norm of each component $T^{(a/q)_d}$ by
\beq\label{T_aq_bound}
 ||T^{(a/q)_d}||_{\ell^p(\Z^k)} \leq A|q_{(d)}|^{-\eta} ,
 \eeq
where 
\beq\label{qs_dfn}
|q_{(d)}| = \sum_{|\al| + |\be| =d} |q_{\al, \be}|,
\eeq
for some $\eta >0$. (In this last sum, $q_{\al, \be}$ ranges over all the denominators in the fixed collection $(a/q)_d$.)
This is sufficient to prove a bound of the form 
\[||T_M||_{\ell^p(\Z^k)} \leq A\]
 for the major operator of level $d$, since as we will see later (Lemma \ref{Obs2_lemma}), the denominators $q$ arising at each level are dyadically separated. 

In order to prove (\ref{T_aq_bound}), we will employ a translation argument that reduces proving $\ell^p(\Z^k)$ bounds for an operator with compactly supported kernel to proving $\ell^p(B)$ bounds for a ``shifted'' operator, where $B$ is a ball with finite radius. 
\begin{lemma}[Shifted ball reduction]\label{shifted_ball_lemma}
Suppose $T$ is an operator acting on functions of $\Z^k$ by
\[ Tf(n) = \sum_{m\in \Z^k} K(n,m) f(m).\]
For any $z \in \Z^k$, define the translated operator $T_z$ by
\[ T_zf(n) =  \sum_{m \in \Z^k} K(n+z,m+z)f(m). \]
Suppose furthermore that the kernel $K(n,m)$ is supported where $|n-m| \leq \rho$, for some fixed radius $\rho$, and let $B_\rho$ denote the ``ball'' of integers $\{n \in \Z^k: |n| \leq \rho\}$. 
Then for any $1 \leq p \leq \infty$, 
\[ ||T||_{\ell^p(B_\tau)} \leq C \sup_{|z| \leq \tau} ||T_z||_{\ell^p(B_\rho)},\]
for any $\rho \leq \tau \leq \infty$, where the constant $C$ depends on the dimension but is independent of the operator $T$.
\end{lemma}
As noted in \cite{SW1}, this lemma is simply a consequence of the fact that the ball $B_\tau$ may be covered by translates of the ball $B_\rho$ with bounded overlap.
In our case, the shifted operator takes the form
\[ T_z^{(a/q)_d}f(n) = \sum_{\bstack{m \in \Z^k}{m-n \neq 0}} \sum_{j \; \text{major}} e^{2\pi i Q(n+z,m+z)}K_j(n-m)f(m),\]
again under the further restrictions on $j$ relating to the fixed collection $(a/q)_d$. Note that the shift by $z$ does not affect the highest degree terms in $Q$. Applying the shifted ball reduction with shifts $z = z_d$,  $\tau = \rho_{d+1} = \infty$, and $\rho = \rho_d$ a finite radius to be chosen below, reduces matters to proving 
\[||T^{(a/q)_d}_{z_d}||_{\ell^p(B_{\rho_d})} \leq A|q_{(d)}|^{-\eta}.\]

We now carry out a major/minor decomposition at level $d-1$ on the shifted operator $T_{z_d}^{(a/q)_d}$, following the same procedure as for level $d$. Once again, the minor components that arise will be bounded directly, while the major component will be further decomposed with respect to collections of fractions $(a/q)_{d-1}$, shifted, and once more subjected to a major/minor decomposition, this time with respect to the coefficients of terms of degree $d-2$ in the polynomial $Q$. Pictorially, this process may be represented as follows, with the major/minor decomposition $M/m$ alternating with the shifting procedure $\tau$: 
\[ T \stackrel{M/m}{\longrightarrow} T^{(a/q)_d} \stackrel{\tau}{\longrightarrow} T^{(a/q)_d}_{z_d} \stackrel{M/m}{\longrightarrow}  T^{(a/q)_{d}, (a/q)_{d-1}}_{z_d} \stackrel{\tau}{\longrightarrow} T^{(a/q)_d, (a/q)_{d-1}}_{z_d,z_{d-1}} \stackrel{M/m}{\longrightarrow}  \cdots \]
At each step the shifts $z_{l}$ lie in a ball $B_{\rho_{l+1}}$ of radius $\rho_{l+1},$ where $\infty=\rho_{d+1} \geq \rho_d \geq \rho_{d-1} \geq \cdots \geq \rho_2.$

This finite inductive procedure stops when all coefficients of terms of degree 2 or greater have been taken into account. The collection of operators that remains at the end of this procedure comprises operators of the form $T^\sharp$:
\beq\label{T_sharp_dfn}
 T^\sharp =  T^{(a/q)_d, (a/q)_{d-1}, \ldots, (a/q)_2}_{z_d,z_{d-1}, \ldots, z_2},
 \eeq
 with a fixed set of collections $(a/q)_d, (a/q)_{d-1}, \ldots, (a/q)_2$ and shifts $z_d,z_{d-1}, \ldots, z_2$.
The key point is that an operator of this form can be factorized, up to acceptable error, as 
\[ T^\sharp = S \otimes T^{\natural} ,\]
where $S$ is a Gauss sum operator and $T^{\natural}$ incorporates the kernel $K$. In Section \ref{sec_induction} we prove that this inductive procedure generalizes to $\ell^p$ results, and in Section \ref{sec_final_bd} we bound the operators $S$ and $T^\natural$.

\section{The inductive procedure}\label{sec_induction}
Having described the inductive procedure in a purely formal manner, we now prove that in order to conclude that
\beq\label{T_p_bound}
||T||_{\ell^p(\Z^k)} \leq A,
\eeq
 it suffices to prove that for each fixed operator $T^\sharp$ of the form (\ref{T_sharp_dfn}),
\beq\label{T_sharp_bound}
 ||T^\sharp||_{\ell^p(B_{\rho_2})} \leq A \prod_{s=2}^d |q_{(s)}|^{-\eta},
 \eeq
for some $\eta >0$ and some finite radius $\rho_2$, where $|q_{(s)}|$ denotes the sum of all denominators of level $s$ in the collection $(a/q)_s$, as defined in (\ref{qs_dfn}).

\subsection{The base case}
We first consider the base case: the passage from $T^\sharp_{d} := T^{(a/q)_d}_{z_d}$ to $T$. This means we must show that (\ref{T_p_bound}) follows from an $\ell^p$ bound for $T_j$ for $j$ minor of level $d$, and the bound
\beq\label{T_sharpd_bound}
|| T^{(a/q)_d}_{z_d}||_{\ell^p(B_{\rho_{d}})} \leq A |q_{(d)}|^{-\eta}
\eeq
for each fixed collection $(a,q)_d$ and shift $z_d \in B_{\rho_{d+1}}$.

\subsubsection{Bound for minor $j$ of level $d$ }\label{T_j_minor}
In order to bound $T_j$ for each minor $j$ of level $d$, i.e. each $j$ for which at least one of the denominators $q_{\al, \be}$ of level $d$ has $q_{\al, \be} > 2^{\ep_d j}$, we will use a Weyl-type bound. We consider exponential sum operators of the form 
\beq\label{S_Om_dfn}
 Sf(n) = \sum_{m \in \Omega} e^{2\pi i P(n,m)} \phi(n,m) f(m).
 \eeq
Here the set $\Omega \subset \Z^k$ is assumed to be of the form $\Omega = \Z^k \cap \omega$ where $\omega $ is a convex set in $\R^k$ contained in a ball of radius $cr$ centered at the origin, for some constant $c$ and parameter $r$. 
Additionally, $P$ is a real-valued polynomial and $\phi$ is a $C^1$ function that satisfies
\beq\label{phi_conditions}
 |\phi(x)| \leq 1, \quad \quad |\nabla \phi(x) | \leq 1/r.
\eeq

\begin{prop}\label{Prop5}
Suppose that for some $\al,\be$ with $|\al| + |\be| \leq d$, $\al \neq 0, \be \neq 0$, the coefficient $\theta_{\al,\be}$ of $P$ has the property that there exist integers $(a_{\al,\be},q_{\al,\be})=1$ and a real number $\ep>0$ such that $|\theta_{\al,\be} - a_{\al,\be}/q_{\al,\be}| \leq 1/q_{\al,\be}^2$, with $r^\ep < q_{\al,\be} \leq r^{|\al|+ |\be| - \ep}.$ Then for every $1< p < \infty$,
\[ ||S||_{\ell^p(\Om)} = O(r^{k - \del}),\]
where $\del = \del(\ep,d,k,p)>0$, but is otherwise independent of the set $\Omega$, the function $\phi$, the fraction $a_{\al,\be}/q_{\al,\be}$, and the coefficients of $P$. 
\end{prop}
When $p=2$, this is Proposition 5 of \cite{SW1}, proved via an $SS^*$ argument and a Weyl-type bound applied to the kernel of $SS^*$. 
When $p=\infty$, note that  trivially
\[ ||Sf(n)||_{\ell^\infty(\Om)} \leq |\Omega| \; ||f||_{\ell^\infty(\Om)} \leq Ar^k ||f||_{\ell^\infty(\Om)},\]
where the constant $A$ is independent of the set $\Omega.$
Thus the full Proposition \ref{Prop5} follows by interpolation and taking adjoints.

Recall that each $T_j$ has kernel supported in a ball of radius $2^{j+1}$. Thus for each minor $j$, we apply Proposition \ref{Prop5} with $r = 2^{j+1}$, $\ep = \ep_d$, and $\phi(n,m) = r^k K_j(n-m)$ to conclude that
\[ ||T_j||_{\ell^p(B_{2^{j+1}})} \leq A2^{-j\del},\]
which by the shifted ball reduction implies that 
\[ ||T_j||_{\ell^p(\Z^k)} \leq A2^{-j\del}.\]
Hence, upon summing over all $j$ that are minor of level $d$, we may conclude that the total minor operator of level $d$ is bounded:
\[ ||T_m||_{\ell^p(\Z^k)} \leq A.\]

\subsubsection{Bound for major $j$ of level $d$}
For the major $j$ we will encounter at each step in the induction procedure, we need another simple lemma, Observation 2 in \cite{SW1}, which states that for a given real number $\theta$, the denominators that occur in Dirichlet approximations to $\theta$ in the major case are dyadically separated.
\begin{lemma}[Dyadic separation]\label{Obs2_lemma}
Given $\theta$, assume that it has two approximations,
\[ |\theta - a/q| \leq 1/qN, \qquad |\theta - a'/q'| \leq 1/q'N',\]
where $1 \leq q \leq N,$ $1 \leq q' \leq N'$. Moreover, suppose that $q \leq N^\ep$, $q' \leq (N')^\ep$ for $\ep$ sufficiently small relative to $N,N'$. Then only one of the three following options occurs: $a/q=a'/q'$, $q \geq 2q'$, or $q' \geq 2q$. 
\end{lemma}
To prove this, we simply note that if $a/q \neq a'/q'$ then
\[ 1/qq' \leq |a/q - a'/q'| \leq 1/qN + 1/q'N',\]
and hence $1 \leq q'/N + q/N'$. One of these summands must be at least $1/2$; supposing $q'/N \geq 1/2$, then $q' \geq N/2 \geq 2N^\ep \geq 2q$, as long as $N \geq 4N^\ep.$ If the other summand is the larger, this leads to $q \geq 2q'$, and the lemma follows. (Note that it is sufficient to take $N,N' >16,$ $\ep<1/2$.) 

We now turn to the major operator of level $d$, defined by $T_M = \sum_{(a/q)_d} T^{(a/q)_d}$. Here we recall that $(a/q)_d = \{a_{\al,\be}/q_{\al,\be} \}$ is a collection of fractions of level $d$ such that for some $j \geq 0$, all the $q_{\al,\be} \leq 2^{\ep_d j}$ and 
\[ |\theta_{\al,\be}  - a_{\al,\be}/q_{\al,\be}| \leq q_{\al,\be}^{-1} 2^{-(d-\ep_d)j}\]
for all $|\al| + |\be|=d$.  We have reduced (\ref{T_p_bound}) to showing  $||T_M||_{\ell^p(\Z^k)} \leq A$, and by Lemma \ref{Obs2_lemma},  the $q_{\al, \be}$ that arise as the collection $(a/q)_d$ varies are dyadically separated, so that it suffices to prove 
\[ ||T^{(a/q)_d}||_{\ell^p(\Z^k)} \leq A |q_{(d)}|^{-\eta},\]
for some $\eta >0$. 
In turn, this follows from (\ref{T_sharpd_bound}) by the shifted ball reduction, with shifts $z_d \in B_{\rho_{d+1}}$ and the radius $\rho_d$ chosen to be\footnote{
 In \cite{SW1}, $\rho_d$ is chosen to be $\rho_d = \inf_{|\al| + |\be| =d} ( |\ga_{\al, \be}|^{-1/(d-\ep_d)})$; however the choice above appears more efficacious in other applications of the method. Both choices satisfy $\rho_d \geq 2^j$.
 } 
\[\rho_d = \inf_{|\al| + |\be| =d} (q_{\al, \be}|\ga_{\al, \be}|)^{-1/(d-\ep_d)},\]
 where $\ga_{\al, \be} = \theta_{\al, \be} - a_{\al, \be}/q_{\al, \be}$.
This completes the proof of the base case.

\subsection{The inductive step}
Next we prove the inductive step: suppose we have reached the operator at the $l$-th stage of the reduction, namely 
\beq\label{Tsharp}
 T_l^\sharp :=  T^{(a/q)_d, (a/q)_{d-1}, \ldots, (a/q)_l}_{z_d,z_{d-1}, \ldots, z_l}.
 \eeq
We must show that the bound
\beq\label{T_sharpl_bound}
 ||T_l^\sharp||_{\ell^p(B_{\rho_l})} \leq A \prod_{s=l}^d |q_{(s)}|^{-\eta},
 \eeq
 follows from the inductive hypothesis that
 \beq\label{T_sharpl'_bound}
   ||T_{l-1}^\sharp||_{\ell^p(B_{\rho_{l-1}})} \leq A \prod_{s=l-1}^d |q_{(s)}|^{-\eta},
   \eeq
for some $\eta >0$, where $T_{l-1}^\sharp$ is defined analogously to (\ref{Tsharp}).

It is worth examining the form the operator $T_l^\sharp$ takes. First, the polynomial $Q$ has been replaced by a polynomial that has been shifted by $z_{s}$ for $l \leq s \leq d$:
\[  \sum_{\bstack{l  \leq s \leq d}{|\al| + |\be| =s}} \theta_{\al, \be} (n+ z_{s})^\al (m + z_{s})^\be + \sum_{|\al| + |\be| < l} \theta_{\al, \be} n^\al m^\be.\]
Here for each $l \leq s \leq d$, the vector $z_{s}$ lies in the ball $B_{\rho_{s+1}}$ of radius $\rho_{s+1}$, where the sequence of radii satisfies $\infty= \rho_{d+1} \geq \rho_d \geq \cdots \geq \rho_l$. (Note: strictly speaking, the coefficients $\theta_{\al,\be}$ that appear in the above sums are no longer the original coefficients of $Q$, but linear combinations of the original coefficients with coefficients depending on shifts at previous stages; however, to simplify notation we continue to use the notation $\theta_{\al,\be}$ at each step.)

The radii $\rho_s$ at the previous stages have been chosen via the Dirichlet approximations of the coefficients. Namely, for each $l \leq s \leq d$ and all $|\al| + |\be|=s$, we have an approximation for each coefficient of degree $s$ as
\[ |\ga_{\al, \be}| = |\theta_{\al, \be}  - a_{\al, \be}/q_{\al, \be} | \leq q_{\al, \be}^{-1} 2^{-(s-\ep_s)j},\]
where $1 \leq q_{\al, \be} \leq 2^{(s-\ep_s)j}$, with $1 \leq a_{\al, \be} \leq q_{\al, \be}$ and $(a_{\al, \be}, q_{\al, \be})=1$.
At each step we have chosen
\beq\label{rho_choice}
\rho_s = \min(\rho_{s+1}, \inf_{|\al| + |\be| = s} (q_{\al, \be} |\ga_{\al, \be}|)^{-1/(s-\ep_s)}).
\eeq

Finally, we note that since $T_l^\sharp$ is the result of a repeated major/minor dichotomy in which at each previous level we preserved only the major $j$, then if we write 
\beq\label{T_lj_sharp}
 T_l^\sharp = \sum_j T^\sharp_{l,j}
 \eeq
 this sum has the additional restrictions on $j$ that for each $l \leq s \leq d$, 
\[q_{\al, \be}|\ga_{\al,\be}| \leq 2^{-(s-\ep_s)j}, \qquad q_{\al, \be} \leq 2^{\ep_sj},\]
or in other words, $2^j \leq c \rho_l$.
Recalling that $j$ corresponds to the dyadic decomposition of the kernel $K = \sum_{j} K_j$, where each $K_j$ is supported where $2^{j-1} \leq |n| \leq 2^{j+1}$, it follows that the kernel of $T_l^\sharp$ is supported where $|n| \leq 2c \rho_l$. 

We are now ready to proceed with decomposing $T_l^\sharp$ via a major/minor dichotomy of level $l-1$. To do so, for each $j$ arising in the sum (\ref{T_lj_sharp}), we find Dirichlet approximations to all coefficients $\theta_{\al,\be}$ in $Q$ with $|\al| + |\be| =l-1$:
\[ |\theta_{\al, \be} - a_{\al, \be}/q_{\al, \be}| \leq q_{\al, \be}^{-1} 2^{-(l-1-\ep_{l-1})j},\]
where $1 \leq q_{\al, \be} \leq 2^{(l-1-\ep_{l-1})j}$ and $1 \leq a_{\al, \be} \leq q$ with $(a_{\al, \be}, q_{\al, \be}) = 1$.
Here $\ep_{l-1}$ is chosen to be small and such that $\ep_l, \ldots, \ep_d$ are small with respect to $\ep_{l-1}$. We now encounter one of two scenarios: either $q_{\al, \be}>2^{(\ep_{l-1})j}$ for some $|\al| + |\be| =l-1$, i.e. $j$ is minor of level $l-1$; or $q_{\al, \be} \leq 2^{(\ep_{l-1})j}$ for all $|\al| + |\be|=l-1$, i.e. $j$ is major of level $l-1$. We will treat each of these cases separately.

\subsubsection{Bound for minor $j$ of level $l-1$}
While Proposition \ref{Prop5} is sufficient for bounding $T_j$ when $j$ is a minor index of the highest level $d$, we need a further variant when bounding $T_j$ for $j$ a minor index at a later stage in the inductive procedure. At these later stages, the polynomial phase has been shifted, and so we consider polynomials of the form
\beq\label{P_shift}
P(n,m) = \sum_{\bstack{|\al| + |\be| \leq d}{|\al| \neq 0, |\be| \neq 0}} \theta_{\al,\be}(n+z_{\al,\be})^\al(m+z_{\al,\be})^\be,
\eeq
for vectors $z_{\al,\be} \in \Z^k$.

\begin{prop}\label{Prop6}
Suppose that for some $\al_0,\be_0$ with $|\al_0| + |\be_0| \leq d$, $\al_0 \neq 0, \be_0 \neq 0$, the coefficient $\theta_{\al_0,\be_0}$ of $P$ has the property that there exist integers $(a_{\al_0,\be_0},q_{\al_0,\be_0})=1$ and a real number $\ep_0>0$ such that $|\theta_{\al_0,\be_0} - a_{\al_0,\be_0}/q_{\al_0,\be_0}| \leq 1/q_{\al_0,\be_0}^2$, with $r^{\ep_0} < q_{\al_0,\be_0} \leq r^{|\al_0|+ |\be_0| - \ep_0}.$ Assume moreover that for each $(\al,\be)$ with $|\al| + |\be|>|\al_0| + |\be_0|$, $\al \neq 0,\be \neq 0$, there exist integers $(a_{\al,\be},q_{\al,\be})=1$ and $\ep>0$ such that $|\theta_{\al,\be} - a_{\al,\be}/q_{\al,\be}| \leq 1/r^{-|\al| - |\be| + \ep}$, with $q_{\al,\be} \leq r^{\ep},$ where $\ep$ is sufficiently small compared to $\ep_0$. Also assume that $|z_{\al,\be}| \leq C|\gamma_{\al,\be}|^{-1/(|\al| + |\be| - \ep)}$, where $\gamma_{\al,\be} = \theta_{\al,\be} - a_{\al,\be}/q_{\al,\be}$.
Then for every $1 < p < \infty$, the operator $S$ defined in (\ref{S_Om_dfn}) with $P$ as in (\ref{P_shift}) is bounded on $\ell^p(\Om)$, with
\[ ||S||_{\ell^p(\Om)} = O(r^{k - \del}),\]
where $\del = \del(\ep_0,\ep, d,k,p)>0$, but is otherwise independent of the set $\Omega$, the function $\phi$, the fraction $a_{\al_0,\be_0}/q_{\al_0,\be_0}$, and the coefficients of $P$. 
\end{prop}
The case $p=2$ is Proposition 6 in \cite{SW1}, which in combination with the trivial $\ell^\infty$ bound yields the $\ell^p$ result for all $1<p< \infty$ by interpolation and taking adjoints.

In order to bound the operator $T_{l,j}^\sharp$ in the case that $j$ is minor of level $l-1$, it is sufficient, by the shifted ball reduction, to prove that 
\beq\label{T_lj_bound'}
 ||T_{l,j}^\sharp||_{\ell^p(B_{\rho_l})} \leq A 2^{-j\del}
 \eeq
for some $\del>0$.  But furthermore, the shifted ball reduction shows that this will be a consequence of proving the bound
\beq\label{T_lj_bound}
  ||T_{l,j}^\sharp||_{\ell^p(\Om +z)} \leq A 2^{-j\del}
  \eeq
uniformly in $z \in B_{\rho_l}$, where $\Om = \{ n: |n| \leq 2^{j+1}\}$. In order to prove (\ref{T_lj_bound}), we apply Proposition \ref{Prop6} with  $|\al_0| + |\be_0| = l-1$, $r=2^{j+1}$, $\phi(n,m) = r^k K_j(n-m)$, and $z_{\al, \be} = z_{s-1} + z$ if $|\al| + |\be| = s$ with $l< s \leq d$, and $z_{\al, \be} = z$ if $|\al| + |\be| =l$.
Having proved (\ref{T_lj_bound}), we then obtain (\ref{T_lj_bound'})
for each minor $j$ of level $l-1$ occurring in (\ref{T_lj_sharp}). But recall that all $j$ appearing in (\ref{T_lj_sharp}) are necessarily major of all previous levels $s$, for all $l \leq s \leq d$, and hence are restricted by the condition $q_{\al, \be} \leq 2^{\ep_s j}$ for all $l \leq s = |\al| + |\be| \leq d$. Therefore in (\ref{T_lj_bound'}), $2^{-j\del} \leq q_{\al, \be}^{-\del'}$ for some small $\del ' >0$ and hence
\[ \sum_{\bstack{j \; \text{minor}}{\text{of level $l-1$}}} ||T_{l,j}^\sharp||_{\ell^p(B_{\rho_l})} 	
	\leq A \prod_{s=l}^d |q_{(s)}|^{-\eta},\]
 for some small $\eta >0$. This is sufficient for the bound (\ref{T_sharpl_bound}). 

\subsubsection{Bound for major $j$ of level $l-1$}
We are thus left with the second scenario: bounding $T_{l,j}^\sharp$ where $j$ is major of level $l-1$, i.e. $q_{\al, \be} \leq 2^{(\ep_{l-1})j}$ for all $|\al| + |\be|=l-1$. We follow the same pattern as we did at level $d$; we will only sketch this step, as the details presented in \cite{SW1} are now identical for both  $\ell^2$ and $\ell^p$ bounds. Let $(a/q)_{l-1}$ denote a collection of rational approximations $\{ a_{\al, \be}/ q_{\al, \be}\}$ to coefficients $\theta_{\al,\be}$ in $Q$ with $|\al| + |\be| = l-1$, and define
\[ T_l^{\sharp, (a/q)_{l-1}}= \sum_j T_{l,j}^\sharp,\]
where $j$ are restricted to those $j$ for which $a_{\al, \be, j}/q_{\al, \be, j} = a_{\al, \be}/q_{\al, \be}$ for all $|\al| + |\be| =l-1$. Then by the dyadic separation of the denominators, the desired bound (\ref{T_sharpl_bound}) will follow from the bound
\beq\label{T_sharp_frxn}
 ||T_l^{\sharp, (a/q)_{l-1}}||_{\ell^p(B_{\rho_l})} \leq A \prod_{s=l-1}^d |q_{(s)}|^{-\eta}.
 \eeq
To reduce this further to the inductive hypothesis (\ref{T_sharpl'_bound}), we need only note that the shifted ball reduction, applied with shifts $z_{l-1} \in B_{\rho_{l}}$,  passes the problem to finding an $\ell^p(B_{\rho_{l-1}})$ norm for $T^{\sharp,(a/q)_{l-1}}_{l,z_{l-1}}$, where $\rho_{l-1}$ is chosen analogously to (\ref{rho_choice}).  This allows us to deduce (\ref{T_sharp_frxn})  from the bound (\ref{T_sharpl'_bound}) for $T_{l-1}^\sharp = T_{l,z_{l-1}}^{\sharp, (a/q)_{l-1}}.$ This completes the inductive step. 

\section{Bounding the final operator}\label{sec_final_bd}
This inductive procedure has reduced the problem to bounding operators of the form $T^\sharp = T^\sharp_2$ defined in (\ref{T_sharp_dfn}), which we will write as $T^\sharp = \sum_j T_j^\sharp$, summed over all the $j$ remaining at this stage. By this final stage, all the coefficients $\theta_{\al, \be}$ corresponding to terms of degrees $2 \leq s \leq d$ have fixed rational approximations, say $\theta_{\al, \be} = a_{\al, \be} /q_{\al, \be} + \ga_{\al, \be}$. Thus the remaining $j$ satisfy for all $2 \leq |s| \leq d$ the conditions
\beq\label{j_restrxns}
q_{\al, \be} |\ga_{\al, \be}| \leq 2^{-(s-\ep)j}, \qquad q_{\al, \be} \leq 2^{\ep_s j},
 \eeq
where in each case $s = |\al| + |\be|$. 
We now choose $Q = \lcm \{ q_{\al, \be} \}$, taken over the finite set of all the denominators corresponding to terms of all degrees. Note that $Q \leq \prod q_{\al, \be} \leq 2^{\ep_0 j}$ for each allowable $j$,
 for some small $\ep_0>0$, by condition (\ref{j_restrxns}) on the denominators.

Our goal is to bound $T^\sharp$ on $\ell^p(B_{\rho_2})$, as in (\ref{T_sharp_bound}). Following \cite{SW1}, $T^\sharp$ may be written, after two simple approximation steps, as:
\[ T^\sharp f(n) = \sum_{l \in (\Z/Q\Z)^k} e^{2\pi i A(r,l)} \sum_{\bar{m} \in B_2^0} e^{2\pi i B(\bar{n},\bar{m})} K(Q(\bar{n}-\bar{m})) f(\bar{m}Q + l) + E_1f(n) + E_2f(n).
 \]
 Here $B_2^0 = \{\bar{n} \in \Z^k: |\bar{n}| \leq (\rho_2 + \sqrt{k})/Q \}$. The integer $k$-tuples $n = \bar{n}Q + r$, $m = \bar{m}Q + l$, which we also denote by $(r,\bar{n})$ and $(l, \bar{m})$, belong to the set $B_2^* = (\Z/Q\Z)^k \times B_2^0$. (Note that $B_{\rho_2} \subseteq B_{2}^*$, and it is ultimately sufficient to prove (\ref{T_sharp_bound}) with $B_{\rho_2}$ replaced by $B_2^*$.)
 The two phases are given by
 \begin{eqnarray}
 A(r,l) &=& \sum_{\bstack{|\al| + |\be| =s}{2 \leq s \leq d}} \frac{a_{\al, \be}}{q_{\al, \be}}( r + z_{s})^\al (l + z_{s})^\be, \label{Arl_dfn}\\
 B(\bar{n}, \bar{m}) & = & \sum_{\bstack{|\al| + |\be| =s}{2 \leq s \leq d}}  \ga_{\al, \be} (\bar{n}Q + z_{s})^\al (\bar{m}Q + z_{s})^\be. \label{Bnm_dfn}
 \end{eqnarray}
 The key point is that $A(r,l)$ is independent of $\bar{n}, \bar{m}$, while $B(\bar{n}, \bar{m})$ is independent of $r,l$; this will allow us to split the main term in $T^\sharp$ into a product $S \otimes T^\natural$.
 
 \subsection{The error term operators}
We first bound the error term operators $E_1, E_2$ that arise in the approximation steps.
As outlined in \cite{SW1}, the first error term operator, acting on functions $f$ of $B_2^*$, takes the form 
\[ |E_1 f(n)| \leq  \sum_j \sum_{2^{j-1} \leq |n-m| \leq 2^{j+1}} \frac{ 2^{-j(1-\ep')}}{(1 + |n-m|)^k}  |f(m)| = \sum_{j} |f| * G_j (n),\]
say, where the $j$ are restricted by (\ref{j_restrxns}) as they are in $T^\sharp$. By Young's inequality, $||f*G_j||_{\ell^p} \leq ||G_j||_{\ell_1} ||f||_{\ell^p}$, so it suffices to note that
\[ ||G_j||_{\ell^1} = \sum_{2^{j-1} \leq |n|  \leq 2^{j+1}} \frac{2^{-j(1-\ep')}}{(1 + |n|)^k}   = O(2^{-j(1-\ep')}).\]
Thus $||E_1||_{\ell^p(B_2^*)}  = O(2^{-j_0(1-\ep')})$, where $j_0$ is the smallest allowable $j$ in $T^\sharp$. But by the restrictions (\ref{j_restrxns}), all allowable $j$ must satisfy $q_{\al, \be} \leq 2^{\ep_s j}$ for all $|\al| + |\be| =s $ with $2 \leq s \leq d$, and hence $O(2^{-j_0(1-\ep')}) = O(q_{\al, \be}^{-(1-\ep')/\ep_s})$ for all denominators of level $s$ with $2 \leq s \leq d$. Therefore
\[ ||E_1||_{\ell^p(B_2^*)} \leq A \prod_{s=2}^d |q_{(s)}|^{-\eta},\]
for some $\eta>0,$ which is of the form (\ref{T_sharp_bound}), as desired.

The second error term operator takes the form
\[ |E_2 f(n)| \leq  \sum_j \sum_{2^{j-1} \leq |n-m| \leq 2^{j+1}} \frac{Q}{(1 + |n-m|)^{k+1}} |f(m)|= \sum_j |f| *H_j,\]
say, where the sum is over allowable $j$ in $T^\sharp$. 
Now 
\[ ||H_j||_{\ell^1} = \sum_{2^{j-1} \leq |n| \leq 2^{j+1}}\frac{Q}{(1 + |n|)^{k+1}}  = O(Q2^{-j}),\]
and recall that $Q \leq 2^{\ep_0 j}$ for some small $\ep_0>0$, so that in fact $O(Q2^{-j}) = O(2^{-(1-\ep_0)j})$.
Thus $ ||E_2||_{\ell^p(B_2^*)}  \leq \sum_j 2^{-(1-\ep_0)j}$, summed over allowable $j$. But again by the restrictions (\ref{j_restrxns}), this shows that
\[ ||E_2||_{\ell^p(B_2^*)} \leq A \prod_{s=2}^d |q_{(s)}|^{-\eta},\]
for some $\eta>0,$ which is also of the desired form (\ref{T_sharp_bound}).
This completes our consideration of the error term operators.

\subsection{Product formulation for $T^\sharp$} 
We are reduced to considering the product operator 
\beq\label{prod_op}
 T^\sharp = S \otimes T^{\natural},
 \eeq
 acting on functions $f(m) = f(\bar{m}Q + l) = f(l,\bar{m})$ belonging to $\ell^p(\Z/Q\Z)^k \otimes \ell^p(B_2^0).$ 
Here $S$ is the Gauss sum operator defined by
\beq\label{Sf_sum}
\ Sf(r) = \frac{1}{Q^k} \sum_{l \in (\Z/Q\Z)^k} e^{2\pi i A(r,l)} f(l),
\eeq
where $A(r,l)$ is as in (\ref{Arl_dfn}). Note that the phase of $S$ involves only the rational approximations to the original coefficients of the polynomial $Q(n,m)$, imbuing $S$ with an arithmetic character.
The second operator $T^{\natural}$ incorporates the Calder\'{o}n-Zygmund kernel, as well as a phase incorporating the error terms in the Dirichlet approximations:
\[ T^{\natural}f(\bar{n}) = \sum_{\bar{m} \in B_2^0} e^{2\pi i B(\bar{n},\bar{m})}Q^kK(Q(\bar{n} - \bar{m}))f(\bar{m}), \]
where $B(\bar{n}, \bar{m})$ is as in (\ref{Bnm_dfn}).
In order to bound $T^\sharp$, it then suffices to bound $S$, $T^\natural$ individually.

 \subsection{Bounding the Gauss sum operator}
The following Weyl-type bound, a direct consequence of Proposition \ref{Prop6}, holds for the operator $S$:
\begin{prop}\label{Prop7}
For $1< p < \infty$, the operator $S$ defined in (\ref{Sf_sum}) satisfies
\[ ||S||_{\ell^p(\Z/Q\Z)^k} \leq CQ^{-\delta},\]
where $\del_p = \del(d,k,p)>0$.
\end{prop}
Recall that $Q = \lcm\{ q_{\al, \be} \} \geq \inf\{q_{\al, \be} \} \geq (\prod_{s=2}^d |q_{(s)}|)^\sig$ for some small $\sig>0$, and hence
\[||S||_{\ell^p(\Z/Q\Z)^k} \leq A \prod_{s=2}^d |q_{(s)}|^{-\eta} \]
for some small $\eta>0$, which is sufficient for (\ref{T_sharp_bound}).

\subsection{Bounding $T^\natural$}
The operator $T^\natural$ is estimated by comparison to its continuous analogue, acting on functions of $\R^k$ by
\[ \mathcal{I}F(x) = \int_{\R^k} e^{2\pi i B(x,y)}Q^k K(Q(x-y))F(y)dy.\]
Given $f \in \ell^p(B_2^0)$, set $f(\bar{n})=0$ for $\bar{n} \not\in B_2^0$ and define a real-variable companion function $F$ on $\R^k$ by setting $F(x) = f(\bar{n})$ for $x$ in the fundamental unit cube in $\R^k$ centered at $\bar{n}$. Then by a simple approximation, for $|x - \bar{n}| \leq 1$,  
\[ |T^\natural f(\bar{n}) -  \mathcal{I} F(x)| \leq A \int_{|x-y|\geq 1} |x-y|^{-k-\ep'}|F(y)|dy,\]
for some small $\ep'>0$.
Since the kernel of the difference operator $T^\natural- \mathcal{I}$ has $L^1(\R^k)$ norm bounded by a constant, 
\[ || T^\natural f - \mathcal{I}F||_{\ell^p(\Z^k)} \leq A ||F||_{L^p(\R^k)} = A||f||_{\ell^p(\Z^k)}.\]
Thus it only remains to bound the action of $\mathcal{I}$ on $L^p(\R^k)$, for all $1< p < \infty$: this is a consequence of work of Ricci and Stein \cite{RSI} on oscillatory singular integrals, since the kernel $Q^kK(Qx)$ satisfies the same Calder\'{o}n-Zygmund conditions as $K$, uniformly in $Q$.\footnote{In fact the original result of \cite{RSI}  is proved for Calder\'{o}n-Zygmund kernels $K$ of critical degree, but the proof may be modified to give the result in the more general case we consider.} This proves that 
\[ ||T^\natural||_{\ell^p(B_2^0)} \leq A,\]
from whence (\ref{T_sharp_bound}) and (\ref{T_p_bound}) follow. This completes the proof of Theorem \ref{SW_gen}.
 
 \subsection*{Acknowledgements}
 The author would like to thank Elias M. Stein for suggesting this area of inquiry and for his generous advice and encouragement. The author was supported in part by the Simonyi Fund and the National Science Foundation, including DMS-0902658 and DMS-0635607, during this research.

%***************************************
\bibliographystyle{mrl}
\bibliography{AnalysisBibliography}
%***************************************

 \end{document}